\RequirePackage{fix-cm}
\documentclass[smallextended]{svjour3}    
\smartqed 
\usepackage{graphicx}

 \journalname{Math. Prog.}
\usepackage[numbers]{natbib}

\usepackage{tikz}
\usepackage{enumerate}

\usepackage{latexsym}
\usepackage{amssymb}
\usepackage{bm}
\usepackage{mathrsfs}
\usepackage{makeidx}
\usepackage{amsmath}
\usepackage{amssymb}
\usepackage{enumitem}

\newcommand{\mb}[1]{\boldsymbol{#1}}

\usepackage{amsfonts}

\def \Q{\mathbb{Q}}
\def \Re{\mathbb{R}}
\def \Z{\mathbb{Z}}

\newtheorem{assumption}{Assumption}

\usepackage{graphicx,multirow,hhline,array}
\usepackage{algorithm}
\usepackage{algpseudocode}

\begin{document}

\title{A Privacy-Aware Distributed Approach for Loosely Coupled Mixed Integer Linear Programming Problems}
\titlerunning{Exact Distributed MILP}

\author{Mohammad Javad Feizollahi}

\institute{
\at Robinson College of Business, Georgia State University, Atlanta, GA 30303, USA.\\
              \email{mfeizollahi@gsu.edu,}      
}

\date{}

\maketitle

\begin{abstract}
In this paper, we propose two exact distributed algorithms to solve mixed integer linear programming (MILP) problems with multiple agents where data privacy is important for the agents. A key challenge is that, because of the non-convex nature of MILPs, classical distributed and decentralized optimization approaches cannot be applied directly to find their optimal solutions. The proposed exact algorithms are based on adding primal cuts and restricting the Lagrangian relaxation of the original MILP problem. We show finite convergence of these algorithms for MILPs with only binary and continuous variables. We test the proposed algorithms on the unit commitment problem and discuss its pros and cons comparing to the central MILP approach. 
\end{abstract}
  


\section{Introduction}\label{Sec:Intro}

Consider the MILP problem 
\begin{equation}\label{eq:MILPAlg}
\begin{split}
z^\text{IP}:= \min\limits_{ \mb{x}_1, \cdots \mb{x}_N}&  \sum_{\nu \in {\cal P} } \mb{c}^\top_\nu \mb{x}_\nu  \\
\text{s.t. } & \mb{x}_\nu \in X_\nu,  \, \forall \nu\in {\cal P}, \\
 & \sum_{\nu\in {\cal P} } \mb{A}_\nu  \mb{x}_\nu  =\mb{b},
\end{split}
\end{equation}
where ${\cal P}=\{1, \cdots, N\}$ is the set of blocks. In reality, there are cases where each of these blocks are governed by a different agent or owner. Each block $\nu$ has its own  $n_\nu$ dimensional vector $\mb{x}_\nu$ of (discrete and continuous) decision variables, and local linear constraints 
\begin{equation}\label{eq:Local-X}
\mb{x}_\nu \in X_\nu,
\end{equation}
where $ X_\nu$ is a linear mixed integer set. Different blocks of the problem \eqref{eq:MILPAlg} are linked to each other via the following linear \emph{coupling constraints}:
\begin{equation}\label{eq:MILPAlg-Link}
\sum_{\nu\in {\cal P} } \mb{A}_\nu  \mb{x}_\nu  =\mb{b}.
\end{equation}
Each $\mb{A}_\nu$  is a $m\times n_\nu$ matrix, for all $ \nu\in {\cal P}$, $\mb{b}$ is a $m$ dimensional vector, where $m$ is the number of coupling constraints \eqref{eq:MILPAlg-Link}. 
If $\mb{A}_\nu$s  are sparse matrices and the number of coupling constraints \eqref{eq:MILPAlg-Link} is relatively small comparing to the total number of local constraints of type \eqref{eq:Local-X}, then we call the problem \eqref{eq:MILPAlg} a \emph{loosely coupled} MILP. In general, relaxing these coupling constraints makes the remaining problem separable and easier.

In the Lagrangian relaxation (LR),  the coupling constraints can be replaced by a linear penalty term in the objective function. Therefore, the LR of MILP \eqref{eq:MILPAlg} will become a separable MILP problem which can be solved in a distributed manner. In contrast to the convex setting, for nonconvex optimization problems such as MILPs, a nonzero duality gap may exist when the coupling constraints are relaxed by using classical Lagrangian dual (LD). In addition to a possible nonzero duality gap, it is not obvious how to obtain optimal Lagrange multipliers and a primal feasible solution by applying LD for MILPs.
 
Augmented Lagrangian dual (ALD) modifies classical LD by appending a nonlinear penalty on the violation of the dualized constraints.  For  MILP \eqref{eq:MILPAlg} under some mild assumptions,  \cite{Feizollahi:2017Augmented} showed asymptotic zero duality gap property of ALD for MILPs when the penalty coefficient is allowed to go to infinity. They also proved that using any norm as the augmenting function with a sufficiently large but finite penalty coefficient closes the duality gap for general MILPs.
 The main drawback of ALD is that the resulting subproblems are not separable because of the nonlinear augmenting functions. To overcome this issue, the alternating direction method of multipliers (ADMM) \citep{Boyd:2011} and related schemes have been developed for convex optimization problems . However, it is not at all clear how to decompose ALD for MILP problems and utilize parallel computation.
Based on ADMM, a heuristic decomposition method was developed in \citep{Feizollahi:2015Large} to solve MILPs arising from electric power network unit commitment problems.

\cite{Bixby:1995} presented a parallel implementation of a branch-and-bound algorithm for mixed 0-1 integer programming problems.
\cite{Ahmed:2013scenario} and \cite{Deng:2017} developed scenario decomposition approaches for 0-1 stochastic programs.
\cite{Munguia:2018} presented a parallel large neighborhood search framework for finding high quality primal solutions for generic MILPs. The approach simultaneously solved a large number of sub-MILPs with the dual objective of reducing infeasibility and optimizing with respect to the original objective.
\cite{Oliveira:2017} proposed a decomposition approach for mixed-integer stochastic programming (SMILP) problems that is inspired by the combination of penalty-based Lagrangian and block Gauss-Seidel methods.

A key challenge is that, because of the non-convex nature of MILPs, classical distributed and decentralized optimization approaches cannot be applied directly to find their optimal solutions. 
In this paper, we propose a distributed approach to solve loosely coupled MILP problems. where each block solves its own modified  LR subproblem iteratively. This approach provides valid lower and upper bounds for the original MILP problem at each iteration. Based on this distributed approach, we develop two exact algorithms which are able to close the gap between lower and upper bounds, and obtain a feasible and optimal solution to the original MILP problem in a finite number of iterations. The proposed exact algorithms are based on adding primal cuts and restricting the Lagrangian relaxation of the original MILP problem. Note that these cuts are not distributable in general.
We test the proposed algorithms on the unit commitment problem and discuss its pros and cons comparing to the central MILP approach. 

This paper is organized as follows. Details of the assumptions and notations are provided in Section \ref{Sec:Prelim}. In Section \ref{Sec:AlgLiterature}, scheme of the  dual decomposition and ADMM as two well known distributed optimization technique are presented. Our distributed MILP approach with two exact algorithms are discussed in Section \ref{Sec:ExactAlgs}.
Experimental results are discussed in Section \ref{Sec:Computations} and conclusions are presented in Section \ref{Sec:Conclusion}.

\section{Preliminaries}\label{Sec:Prelim}
Let $\Re$, $\Z$, and $\Q$ denote the sets of real, integer and rational numbers, respectively. For a finite dimensional vector $\mb{a}$, denote its transpose by $\mb{a}^\top$. For a set ${\cal S}$, denote its cardinality by $|{\cal S}|$. 
In this paper, we consider MILP problem \eqref{eq:MILPAlg} which satisfies the following assumptions.

\begin{assumption}\label{Assump:MILPAlg}
For the MILP  \eqref{eq:MILPAlg} we have the following:
\begin{enumerate}[label=(\alph*)]

\item For each block $\nu\in {\cal P}$, $X_\nu$ is a linear mixed integer set defined by 
\begin{equation}\label{eq:DefineX}
X_\nu:=\{ (\mb{u}_\nu^\top,\mb{y}_\nu^\top)^\top : \mb{u}_\nu \in U_\nu, ~ \mb{y}_\nu\in Y_\nu(\mb{u}_\nu) \},
\end{equation}
where $\mb{u}_\nu\in \{0,1\}^{n_\nu^1}$ and $\mb{y}_\nu \in \Re^{n_\nu^2}$ are the subvectors of $n_\nu^1$ binary and $n_\nu^2$ continuous decision variables, respectively, with $n_\nu=n_\nu^1+n_\nu^2$.

\item In description \eqref{eq:DefineX} of $X_\nu$, $U_\nu$ and $Y_\nu(\mb{u}_\nu)$ are subsets of  $\{0,1\}^{n_\nu^1}$ and $\Re^{n_\nu^2}$, respectively. Because $U_\nu$ is a finite set, it can be represented by a set of linear inequalities and integrality constraints. For a given $\mb{u}_\nu\in U_\nu$, we assume $Y_\nu(\mb{u}_\nu)$ is a (possibly empty) polyhedron. 
In particular, let $Y_\nu(\mb{u}_\nu)=\{ \mb{y}_\nu: \Re^{n_\nu^2}: \mb{E}_\nu \mb{u}_\nu +\mb{F}_\nu \mb{y}_\nu \le \mb{g}_\nu\}$, where $\mb{E}_\nu $ and $\mb{F}_\nu$ are matrices and $ \mb{g}_\nu$ is a vector of appropriate finite dimensions, independent of the value of $\mb{u}_\nu$. 

\item  $\mb{c}_\nu$, $\mb{A}_\nu$, $\mb{E}_\nu $, $\mb{F}_\nu$ and $ \mb{g}_\nu$, for all $\nu\in {\cal P}$, and $\mb{b}$ have rational entries. 

\item Problem \eqref{eq:MILPAlg} is feasible and its optimal value is bounded.
\end{enumerate}
\end{assumption}

Let $n^1:=\sum_{\nu \in {\cal P} } n_\nu^1$ and $n^2:=\sum_{\nu \in {\cal P} } n_\nu^2$ denote total number of binary and continuous variables, respectively, and $n=n^1+n^2$. For convenience, let 
\begin{equation}\nonumber
\begin{split}
 & \mb{c}:= \left[ \begin{array}{c}
\mb{c}_1 \\ 
 \vdots \\ 
 \mb{c}_N
 \end{array}  \right],~
  \mb{x}  := \left[ \begin{array}{c}
\mb{x}_1 \\ 
 \vdots \\ 
 \mb{x}_N
 \end{array}  \right],~ 
  \mb{u}:= \left[ \begin{array}{c}
\mb{u}_1 \\ 
 \vdots \\ 
 \mb{u}_N
 \end{array}  \right], ~
  \mb{y}:= \left[ \begin{array}{c}
\mb{y}_1 \\ 
 \vdots \\ 
 \mb{y}_N
 \end{array}  \right],\\ 
& \mb{A}:=  [\mb{A}_1,\cdots, \mb{A}_N],~
 X:= X_1\times \cdots \times X_N,~
 U:= U_1\times \cdots \times U_N,\\
 & Y(\mb{u}):= Y_1(\mb{u}_1)\times \cdots Y_N(\mb{u}_N).
\end{split}
\end{equation}
Then, problem \eqref{eq:MILPAlg} can be recast as $z^\text{IP}= \min\limits_{ \mb{x} } \{ \mb{c}^\top \mb{x} : \mb{x}\in X,  \mb{A}  \mb{x} =\mb{b}\}$.

By Assumption \ref{Assump:MILPAlg}-d), there exists a solution $\mb{x}^\ast$ which satisfies constraints \eqref{eq:Local-X} and \eqref{eq:MILPAlg-Link}, and $\mb{c}^\top \mb{x}^\ast =z^\text{IP}$. 
Therefore, by data rationality assumption in part (c), the value of the linear programming (LP) relaxation ($z^\text{LP}$) of \eqref{eq:MILPAlg} is bounded \citep{Blair:1979}, i.e. $-\infty < z^\text{LP} \le z^\text{IP}< \infty$.

\begin{example} \label{Ex:MILPAlgExample}
Following is an example for problem \eqref{eq:MILPAlg} with two blocks.
\begin{equation}\label{eq:MILPAlgExample}
\begin{split}
\min ~~& 70u_{11}+70u_{12}+110u_{13}+2y_{11}+2y_{12}+48u_{21}+48u_{22}+52u_{23}+3y_{21}+3y_{22}\\
\text{s.t. } & 
\left. \begin{array}{c}
u_{12}-u_{11} - u_{13}\le 0,\\
 30u_{11} \le y_{11} \le 100 u_{11},\\
 30u_{12} \le y_{12} \le 100 u_{12},\\
-35 \le y_{12}-y_{11} \le 35,\\
u_{11}, u_{12}, u_{13}\in \{0,1\},
\end{array} \right\} \text{Local constraints for block 1}\\
& \left. \begin{array}{c}
 u_{22}-u_{21} - u_{23}\le 0,\\
 20u_{11} \le y_{21} \le 80 u_{21},\\
20u_{12} \le y_{22} \le 80 u_{22},\\
-30 \le y_{22}-y_{21} \le 30,\\
u_{21}, u_{22}, u_{23}\in \{0,1\},\\
\end{array} \right\} \text{Local constraints for block 2}\\
& \left. \begin{array}{c}
y_{11}+y_{21}= 90,\\
y_{12}+y_{22}= 120.
\end{array}~~~~~~~~~ \right\} \text{Coupling constraints}
\end{split}
\end{equation}
Recalling the notations described in Sections \ref{Sec:Intro} and \ref{Sec:Prelim},
$\mb{u}_1=(u_{11}, u_{12}, u_{13})^\top$ and  $\mb{u}_2=(u_{21}, u_{22}, u_{23})^\top$ are the vectors of binary variables for blocks 1 and 2, respectively. Similarly, $\mb{y}_1=(y_{11}, y_{12})^\top$ and  $\mb{y}_2=(y_{21}, y_{22})^\top$ are the vectors of continuous variables for blocks 1 and 2, respectively.
Then, $\mb{x}_1=(\mb{u}_1^\top,\mb{y}_1^\top)^\top$ and $\mb{x}_2=(\mb{u}_2^\top,\mb{y}_2^\top)^\top$ are the vectors of decision variables for blocks 1 and 2, respectively.  
Moreover, $\mb{u}=(\mb{u}_1^\top,\mb{u}_2^\top)=(u_{11}, u_{12}, u_{13},u_{21}, u_{22}, u_{23})^\top$ and $\mb{y}=(\mb{y}_1^\top,\mb{y}_2^\top)=(y_{11}, y_{12},y_{21}, y_{22})^\top$ are the overall vectors of binary and continuous variables. In this example, we have
$$ \mb{c}_1=\left[
\begin{array}{c}
70\\ 70\\ 110\\ 2\\ 2
\end{array} \right],~
\mb{c}_2=
\left[ \begin{array}{c}
48 \\ 48 \\ 52 \\ 3 \\ 3
\end{array} \right], \text { and } 
 \mb{A}_1= \mb{A}_2=\left[
\begin{array}{c c c c c}
0 & 0 & 0 & 1 & 0\\
0 & 0 & 0 & 0 & 1
\end{array} \right]. 
$$
Moreover,
$$U_1=\{\mb{u}_1\in \{0,1\}^3: u_{12}-u_{11} - u_{13}\le 0 \},$$
$$U_2=\{\mb{u}_2\in \{0,1\}^3: u_{22}-u_{21} - u_{23}\le 0 \},$$
$$U=U_1\times U_2=\left\{\mb{u}\in \{0,1\}^6: \begin{array}{c} u_{12}-u_{11} - u_{13}\le 0,\\ u_{22}-u_{21} - u_{23}\le 0 \end{array}\right\},$$
$$Y_1(\mb{u}_1)=\left\{\mb{y}_1\in \Re^2: 
\begin{array}{c}
 30u_{11} \le y_{11} \le 100 u_{11},\\
 30u_{12} \le y_{12} \le 100 u_{12},\\
-35 \le y_{12}-y_{11} \le 35
\end{array}\right \},$$
$$Y_2(\mb{u}_2)=\left\{\mb{y}_2\in \Re^2: 
\begin{array}{c}
 20u_{11} \le y_{21} \le 80 u_{21}\\
20u_{12} \le y_{22} \le 80 u_{22}\\
-30 \le y_{22}-y_{21} \le 30
\end{array}\right \},$$
$$X_1=\left\{\mb{x}_1=(\mb{u}_1^\top,\mb{y}_1^\top)^\top\in \{0,1\}^3\times \Re^2: 
\begin{array}{c}
u_{12}-u_{11} - u_{13}\le 0,\\
 30u_{11} \le y_{11} \le 100 u_{11},\\
 30u_{12} \le y_{12} \le 100 u_{12},\\
-35 \le y_{12}-y_{11} \le 35,
\end{array}\right \},$$
$$X_2=\left\{\mb{x}_2=(\mb{u}_2^\top,\mb{y}_2^\top)^\top\in \{0,1\}^3\times \Re^2: 
\begin{array}{c}
 u_{22}-u_{21} - u_{23}\le 0,\\
 20u_{11} \le y_{21} \le 80 u_{21},\\
20u_{12} \le y_{22} \le 80 u_{22},\\
-30 \le y_{22}-y_{21} \le 30,
\end{array}\right \}.$$
\end{example}

For a given vector of the dual (Lagrange) variables, $\mb{\mu}\in \Re^m$, the standard LR for MILP \eqref{eq:MILPAlg} is
\begin{equation}\label{eq:MILP1-LR}
\begin{split}
z^\text{LR}(\mb{\mu}):=\mb{\mu}^\top \mb{b}+ \min\limits_{ \mb{x}_1,\cdots, \mb{x}_N }  & \sum_{ \nu\in {\cal P} } {\cal L}_\nu(\mb{x}_\nu,\mb{\mu})
\\
\text{s.t. } &\mb{x}_\nu \in X_\nu,  \, \forall \nu\in {\cal P}, 
\end{split}
\end{equation}
where
\begin{equation}\nonumber
{\cal L}_\nu(\mb{x}_\nu,\mb{\mu}):=(\mb{c}^\top_\nu -\mb{\mu}^\top \mb{A}_\nu ) \mb{x}_\nu,  \, \forall \nu\in {\cal P}, 
\end{equation}
and the corresponding LD value is 
\begin{equation}\label{eq:MILP1-LD}
z^\text{LD}:= \sup\limits_{\mb{\mu}\in \Re^m} z^\text{LR}(\mb{\mu}).
\end{equation}
Since \eqref{eq:MILP1-LR} is a relaxation of \eqref{eq:MILPAlg}, $z^\text{LR}(\mb{\mu}) \le z^\text{LD} \le z^\text{IP}$ holds, for any $\mb{\mu}\in \Re^m$. Due to the presence of binary variables, a nonzero duality gap may exists \citep{Wolsey:1999}, i.e. $z^\text{LD} < z^\text{IP}$ is possible. Let $\mb{\mu}^\ast$ be a maximizer in \eqref{eq:MILP1-LD}, which  exists under Assumption \ref{Assump:MILPAlg}.
Obtaining $\mb{\mu}^\ast$ and $z^{LD}$ are not straight forward in practice. A popular and easy approach to solve \eqref{eq:MILP1-LD} is the subgradient decent method, where the problem \eqref{eq:MILP1-LR} is solved iteratively and the dual multipliers are updated at each iteration. Note that problem \eqref{eq:MILP1-LR} is separable and it can be solved by computing 
$$\min\limits_{\mb{x}_\nu}  \{  {\cal L}_\nu(\mb{x}_\nu,\mb{\mu}): \mb{x}_\nu \in X_\nu\}$$
for each block $\nu$.
Even with $\mb{\mu}^\ast$ at hand, a {\em primal feasible solution}, one that satisfies all constraints in model \eqref{eq:MILPAlg}, is not readily available. 
In other words, an optimal solution of LR \eqref{eq:MILP1-LR} for $\mb{\mu}^\ast$
does not necessarily satisfy the coupling constraints \eqref{eq:MILPAlg-Link} in problem \eqref{eq:MILPAlg}.

For a given $\hat{\mb{u}} \in U$, the best corresponding primal feasible solution, if there exists one, and its objective value, $z(\mb{\hat u})$, can be computed by solving the following LP:
\begin{equation}\label{eq:LP1}
\begin{split}
z(\mb{\hat u}):=  \min\limits_{ \mb{y}_1,\cdots,\mb{y}_N } & ~\sum_{\nu \in {\cal P}}  \mb{c}_\nu^\top \left[ \def\arraystretch{0.5} \begin{array}{c}
\mb{\hat u}_\nu \\ 
\mb{y}_\nu
\end{array} 
\right] \\
\text{s.t. } & \mb{y}_\nu\in Y_\nu(\mb{\hat u}_\nu), \, \forall \nu \in {\cal P},\\
& \sum_{\nu \in {\cal P}}  \mb{A}_\nu \left[ \def\arraystretch{0.5} \begin{array}{l}
\mb{\hat u}_\nu \\ 
\mb{y}_\nu
\end{array} 
\right] =\mb{b}  .
\end{split}
\end{equation}
Problem \eqref{eq:LP1} is an LP and can be solved with a distributed algorithm \citep{Boyd:2011}. 
Denote the upper and lower bounds on $z^\text{IP}$ by $ub$ and $lb$, respectively. Then, $z(\mb{u})$ and  $z^\text{LR}(\mb{\mu})$ are valid $ub$ and $lb$, respectively, for all $\mb{u}\in U$ and $\mb{\mu}\in \Re^m$, i.e.
$$z^\text{LR}(\mb{\mu}) \le z^\text{IP} \le z(\mb{u}),~~ \forall \mb{u}\in U, \mb{\mu}\in \Re^m.$$ 
In fact, 
\begin{equation}\label{eq:z-U}
z^\text{IP}= \min\limits_{\mb{u}\in U} ~  z(\mb{u}).
\end{equation}

\section{Dual Decomposition and ADMM for MILPs} \label{Sec:AlgLiterature}
Dual decomposition and ADMM are two well known distributed optimization technique in the context of convex optimization. Our distributed MILP algorithms in this paper are based on extensions of these two techniques. Next, we present these schemes and discuss challenges in applying them to MILPs.

\subsection{Dual Decomposition}
Dual decomposition is a well known technique to solve large scale optimization problems. Early works on application of dual decomposition for 
large scale linear programming can be found in \citep{Benders:1962,Dantzig:1960,Dantzig:1963,Everett:1963}.

Let $\rho_\mu^k>0$ be the step size for updating the dual vector $\mb{\mu}$ at iteration $k$. Algorithm \ref{Alg:BasicDualD} represents an overall scheme of a dual decomposition method to solve \eqref{eq:MILPAlg}. Each iteration of this method requires a ``broadcast'' and a ``gather'' operation. Dual update step (line \ref{Line:DualUpdate} in Algorithm \ref{Alg:BasicDualD}) requires $\mb{A}_\nu  \mb{x}_\nu^k$ values from all blocks. Once $\mb{\mu}^{k}$ is computed, it must be broadcasted to all blocks.

A lower bound for $z^\text{IP}$ can be obtained from Algorithm \ref{Alg:BasicDualD}. If $\sum\limits_{\nu\in {\cal P} } \mb{A}_\nu  \mb{x}_\nu^k =\mb{b}$ in some iteration $k$ of  this algorithm, $\mb{x}^k$ is a feasible and optimal solution of \eqref{eq:MILPAlg}. But, this case is not likely in practice and there is no hope to find a feasible solution for \eqref{eq:MILPAlg} by running only Algorithm \ref{Alg:BasicDualD}. Therefore, in general we cannot expect to get an upper bound for $z^\text{IP}$ from this algorithm. A modified version of dual decomposition technique is presented in Algorithm \ref{Alg:DualD} which is able to provide upper bounds for $z^\text{IP}$.

\begin{algorithm}
\caption{Basic Dual Decomposition}\label{Alg:BasicDualD}
\begin{algorithmic}[1]
\State $lb \leftarrow -\infty$, $\mb{\mu}^0 \leftarrow \mb{0}$, and $k \leftarrow 0$.
\While{ some termination criteria is not met}
\State $k \leftarrow k+1$
\For{$\nu:=1$ to $N$}
\State solve $\min\limits_{\mb{x}_\nu}  \{  {\cal L}_\nu(\mb{x}_\nu,\mb{\mu}^{k-1}): \mb{x}_\nu \in X_\nu\}$
\State let $v_\nu^k$ be the optimal value and $\mb{x}_\nu^k$ be an optimal solution
\EndFor
\If{$lb<  \mb{\mu}^\top \mb{b}+\sum\limits_{\nu\in {\cal P}} v_\nu^k$}
\State $lb\leftarrow  \mb{\mu}^\top \mb{b}+\sum\limits_{\nu\in {\cal P}}  v_\nu^k$
\EndIf
\State  $\mb{\mu}^{k} \leftarrow \mb{\mu}^{k-1}+\rho_\mu^{k}\left(\mb{b}-
\sum\limits_{\nu\in {\cal P} } \mb{A}_\nu  \mb{x}_\nu^k  \right) $ \label{Line:DualUpdate}
\EndWhile
\end{algorithmic}
\end{algorithm}

\subsection{Alternating Direction Method of Multipliers (ADMM)}\label{Sec:ADMM}
ADMM is an algorithm that is intended to blend the separability of dual decomposition with the superior convergence properties of the method of multipliers \citep{Boyd:2011}. For $\rho>0$ and ${\mb \mu}\in \Re^m$, the augmented Lagrangian with squared Euclidean norm has the following form.
\begin{equation} \label{eq:ALMILP}
{\cal L}_\rho^+(\mb{x}_1,\cdots,\mb{x}_N,{\mb \mu})=  \sum_{\nu\in {\cal P} }  \mb{c}^\top_\nu \mb{x}_\nu  + {\mb \mu}^\top \left(\mb{b}-
\sum_{\nu\in {\cal P} } \mb{A}_\nu  \mb{x}_\nu  \right)
 +\frac{\rho}{2}\left\|\mb{b} -\sum_{\nu\in {\cal P} } \mb{A}_\nu  \mb{x}_\nu  \right\|_2^2.
\end{equation}

A robust relaxation for MILP \eqref{eq:MILPAlg} is the augmented Lagrangian relaxation (ALR) which has the following form:
\begin{equation} \label{eq:ALR} 
\begin{split}
z^\text{LR+}_\rho({\mb \mu}):= \min\limits_{ \mb{x}_1, \cdots \mb{x}_N}&  {\cal L}_\rho^+(\mb{x}_1,\cdots,\mb{x}_N,{\mb \mu}) \\
\text{s.t. } & \mb{x}_\nu \in X_\nu,  \, \forall \nu\in {\cal P},
\end{split}
\end{equation}
and the corresponding ALD value is 
\begin{equation}\label{eq:MILP1-ALD}
z^\text{LD+}:= \sup\limits_{\mb{\mu}\in \Re^m} z^\text{LR+}(\mb{\mu}).
\end{equation}
Since \eqref{eq:ALR} is a relaxation of \eqref{eq:MILPAlg}, $z^\text{LR+}(\mb{\mu}) \le z^\text{LD+} \le z^\text{IP}$ holds, for any $\mb{\mu}\in \Re^m$.

For  MILP \eqref{eq:MILPAlg} under Assumption \ref{Assump:MILPAlg}, \cite{Feizollahi:2017Augmented} showed that
using ALD with any norm as the augmenting function is able to
close the duality gap with a finite penalty coefficient $\rho$.
It is obvious that ${\cal L}_\rho^+ (\mb{x}_1,\cdots,\mb{x}_N,{\mb \mu})$ in \eqref{eq:ALMILP} is not separable between different blocks, because the nonlinear (quadratic) terms are coupling different block to each other. For convex optimization problems, a decomposable algorithm to solve \eqref{eq:ALR} is ADMM \citep{Boyd:2011}.

\subsubsection{ADMM with two blocks}
Algorithm \ref{Alg:ADMM2} presents an ADMM approach for an optimization problem with two blocks. In $k$th iteration of this algorithm, ${\cal L}_\rho^+ (\mb{x}_1, \mb{x}_2^{k-1},\mb{\mu}^{k-1})$ is first minimized with respect to $\mb{x}_1$, assuming that $\mb{x}_2$ is fixed at its previous value $\mb{x}_2^{k-1}$. Then, ${\cal L}_\rho^+ ( \mb{x}_1^{k}, \mb{x}_2,\mb{\mu}^{k-1})$ is minimized with respected to $\mb{x}_2$, assuming that $\mb{x}_1$ is fixed at its previous value $\mb{x}_1^{k}$. Finally, the vector of dual variables $\mb{\mu}^{k}$ is updated. Note that $\rho>0$ is a given and fixed penalty factor.

\begin{algorithm}
\caption{ADMM procedure for two blocks}\label{Alg:ADMM2}
\begin{algorithmic}[1]
\State $\mb{x}_2^0 \leftarrow \mb{0} $, $\mb{\mu}^0 \leftarrow \mb{0} $, and $k\leftarrow 0$
\While{some termination criteria is not met }
\State $k\leftarrow k+1$
\State $\mb{x}_1^{k} \leftarrow \arg\min\limits_{\mb{x}_1 \in X_1} {\cal L}_\rho^+ (\mb{x}_1, \mb{x}_2^{k-1},\mb{\mu}^{k-1})$
\State $\mb{x}_2^{k} \leftarrow \arg\min\limits_{\mb{x}_2 \in X_2} {\cal L}_\rho^+ ( \mb{x}_1^{k}, \mb{x}_2,\mb{\mu}^{k-1})$
\State Update $\mb{\mu}^{k} \leftarrow  \mb{\mu}^{k-1}+ \rho\times  [\mb{b}- (\mb{A}_1  \mb{x}_1^{k}+ \mb{A}_2 \mb{x}_2^{k})]$
\EndWhile
\end{algorithmic}
\end{algorithm}

Let $\mb{\alpha}^k$ and $\mb{\beta}^k$ denote vectors of primal and dual residuals at iteration $k$. Then,
\begin{equation}\nonumber
\begin{split}
\mb{\alpha}^k=\mb{b}- (\mb{A}_1  \mb{x}_1^{k}+ \mb{A}_2 \mb{x}_2^{k}) \text { and }
\mb{\beta}^k=\rho \mb{A}_1^\top \mb{A}_2 (\mb{x}_2^k-\mb{x}_2^{k-1}).
\end{split}
\end{equation}

If problem \eqref{eq:MILPAlg} is solvable and the sets $X_1$ and $X_2$ are convex, closed, and non-empty, Algorithm \ref{Alg:ADMM2} can solve \eqref{eq:MILPAlg} in a distributed framework \citep{Boyd:2011}. In this case, primal residuals ($\mb{\alpha}^k$) converge to zero. Moreover, dual variables ($\mb{\mu}^{k}$) and objective value converge to their optimal values \citep{Boyd:2011}. Note that   discrete variables destroy the nice convergence properties of ADMM for MILP problems \citep{Feizollahi:2015Large}.  
In practice, ADMM converges to modest accuracy --sufficient for many applications-- within a few tens of iterations \citep{Boyd:2011}. However, direct extension of ADMM for multi-block convex
minimization problems is not necessarily convergent \citep{Chen:2016direct}.

\subsubsection{Global Variable Consensus Problem with ADMM}
To extend ADMM for multi-block minimization problems, a global variable consensus problem can be constructed. An equivalent optimization problem for \eqref{eq:MILPAlg} is as follows.
\begin{subequations}\label{eq:consensus}
\begin{eqnarray}
z^\text{IP}:=  \min\limits_{ \mb{x}_1, \cdots \mb{x}_N, \bar{\mb x}_1, \cdots, \bar{\mb x}_N} &&  \sum_{\nu \in {\cal P} } \mb{c}^\top_\nu \mb{x}_\nu \nonumber \\
\text{s.t. } & & \mb{x}_\nu \in X_\nu ,  \, \forall \nu\in {\cal P}, \nonumber \\
& & \sum_{\nu\in {\cal P} } \mb{A}_\nu  \bar{\mb x}_\nu  =\mb{b}, \label{eq:LinkCons1}\\
&& \bar{\mb x}_\nu=\mb{x}_\nu,  \, \forall \nu\in {\cal P}.\label{eq:LinkCons2}
\end{eqnarray}
\end{subequations}

\begin{algorithm}
\caption{Consensus ADMM} \label{Alg:ADMM}
\begin{algorithmic}[1]
\State $\bar{\mb x}^0 \leftarrow \mb{0}$, $\mb{\mu}^0 \leftarrow \mb{0} $, and $k\leftarrow 0$
\While{some termination criteria is not met }
\State  $k\leftarrow k+1$
\For{$\nu:=1$ to $N$}
\State $\mb{x}_\nu^{k} \leftarrow \arg\min\limits_{\mb{x}_\nu \in X_\nu} {\cal L}_{\rho,\nu}^+(\mb{x}_\nu,\bar{\mb x}_\nu^{k-1},\mb{\mu}_\nu^{k-1})$
\EndFor
\State $\bar{\mb x}^{k} \leftarrow \arg\min\limits_{\bar{\mb x}}  \left\{  {\cal L}_{\rho}^+(\mb{x}^{k},\bar{\mb x},\mb{\mu}^{k-1}):  \sum_{\nu\in {\cal P} } \mb{A}_\nu  \bar{\mb x}_\nu  =\mb{b} \right\}$ by using \eqref{eq:Update-xbar}
\For{$\nu:=1$ to $N$}
\State  $\mb{\mu}_\nu^{k} \leftarrow  \mb{\mu}_\nu^{k-1}+ \rho\times  (\mb{x}_\nu^{k}-\bar{\mb x}_\nu^{k})$
\EndFor
\EndWhile
\end{algorithmic}
\end{algorithm}

Formulation \eqref{eq:consensus} can be decomposed into two parts, where one part includes variable vectors $\mb{x}_1, \cdots \mb{x}_N$, constraints $\mb{x}_\nu \in X_\nu $, for all $\nu\in {\cal P}$ and the objective function, and the other part contains variable vectors $\bar{\mb x}_1, \cdots, \bar{\mb x}_N$ and constraints \eqref{eq:LinkCons1}. In this case, constraints \eqref{eq:LinkCons2} are coupling these two parts and  Algorithm \ref{Alg:ADMM2}, ADMM with two blocks, can be adjusted to solve problem \eqref{eq:consensus} in a distributed manner. Algorithm \ref{Alg:ADMM}, consensus ADMM, represents this process.

Let 
\begin{equation}\nonumber
{\cal L}_{\rho}^+(\mb{x},\bar{\mb x},\mb{\mu}):=\sum_{\nu \in {\cal P} } {\cal L}_{\rho,\nu}^+(\mb{x}_\nu,\bar{\mb x}_\nu,\mb{\mu}_\nu),
\end{equation}
where 
$
{\cal L}_{\rho,\nu}^+(\mb{x}_\nu,\bar{\mb x}_\nu,\mb{\mu}_\nu):=\mb{c}^\top_\nu \mb{x}_\nu   +\mb{\mu}_\nu^\top(\mb{x}_\nu-\bar{\mb x}_\nu)
+ \frac{\rho}{2} \|\mb{x}_\nu-\bar{\mb x}_\nu\|_2^2.
$
Then,  the subproblem for part one is 
$
\min\limits_{ \mb{x}} \{ {\cal L}_{\rho}^+(\mb{x},\bar{\mb x},\mb{\mu}) :\mb{x}_\nu \in X_\nu ,  \, \forall \nu\in {\cal P} \},
$
which is separable between blocks and can be solved in parallel. Moreover, the subproblem for part two is 
\begin{equation}\nonumber
\min\limits_{ \bar{\mb x}} \{ {\cal L}_{\rho}^+(\mb{x},\bar{\mb x},\mb{\mu}) :\sum_{\nu\in {\cal P} } \mb{A}_\nu  \bar{\mb x}_\nu  =\mb{b} \}
\end{equation}
which has a closed form solution as follows (assuming $\mb{A}$ has full row rank):
\begin{equation}\label{eq:Update-xbar}
\begin{split}
\arg\min\limits_{\bar{\mb x}} & \left\{ {\cal L}_{\rho}^+(\mb{x},\bar{\mb x},\mb{\mu}):  \sum_{\nu\in {\cal P} } \mb{A}_\nu  \bar{\mb x}_\nu  =\mb{b} \right\}
  =\arg\min\limits_{\bar{\mb x}} \left\{\|\mb{x}+\frac{\mb{\mu}}{\rho}-\bar{\mb x}\|_2^2  :  \mb{A}\bar{\mb x}={\mb b} \right\}\\
&= [I-\mb{A}^\top (\mb{A}\mb{A}^\top)^{-1}\mb{A}](\mb{x}+\frac{\mb{\mu}}{\rho})+\mb{A}^\top (\mb{A} \mb{A}^\top)^{-1}\mb{b}
\end{split}
\end{equation}
where the second equality is well known in linear algebra for finding the orthogonal projection of a point onto an affine subspace  \citep[e.g.][]{Meyer:2000,Plesnik:2007finding}. In general, to compute inverse matrices is not easy \citep{Higham:2002}, but it can be done efficiently for sparse matrices with specific structures.

In distributed consensus optimization, ADMM has a linear convergence rate \citep{Shi:2014linear}. Consensus ADMM can be interpreted as a method for solving problems in which the objective and constraints are distributed across multiple processors. Each processor only has to handle its own objective and constraint term, plus a quadratic term which is updated each iteration. The linear parts of the quadratic terms are updated in such a way that the variables converge to a common value, which is the solution of the full problem \citep{Boyd:2011}.

In our context of MILP \eqref{eq:MILPAlg}, consensus ADMM (Algorithm \ref{Alg:ADMM}) can be used for upper bounding $z^{IP}$. For a given set $\hat{{\cal S}} \subset U$, an upper bounding method is as Algorithm \ref{Alg:UpperBounding}.
\begin{algorithm}
\caption{Upper Bounding Algorithm}\label{Alg:UpperBounding}
\begin{algorithmic}[1]
\For{$\hat{\mb{u}}\in \hat{{\cal S}}$}
\State compute $z(\mb{\hat u})$ by solving LP \eqref{eq:LP1} with consensus ADMM, Algorithm \ref{Alg:ADMM}
\If{$z(\mb{\hat u})<ub$ }
\State $ub \leftarrow z(\mb{\hat u})$
\State $\mb{u}^\ast \leftarrow \mb{\hat u}$
\EndIf
\EndFor
\end{algorithmic}
\end{algorithm}

\subsection{Combination of Dual Decomposition and Consensus ADMM}
A combination of Algorithm \ref{Alg:BasicDualD} (dual decomposition) and Algorithm \ref{Alg:ADMM} (consensus ADMM) can be used to generate lower and upper bounds for $z^\text{IP}$. Algorithm \ref{Alg:DualD} presents a modified version of Algorithm \ref{Alg:BasicDualD}. In this algorithm, for a given binary vector $\hat{\mb u}$, Algorithm \ref{Alg:ADMM} (consensus ADMM) is used to refine continuous variables $\mb{y}$, and obtain an upper bound for $z^\text{IP}$.
Besides the issues related to the non zero duality gap and the challenges in finding the the best dual vector $\mb{\mu}^\ast$, which is a maximizer in \eqref{eq:MILP1-LD}, it is possible for Algorithms \ref{Alg:BasicDualD} and \ref{Alg:DualD} to cycle between non-optimal solutions forever.

\begin{algorithm}
\caption{Modified Dual Decomposition for MILPs}\label{Alg:DualD}
\begin{algorithmic}[1]
\State $ub \leftarrow +\infty$, ${\cal S} \leftarrow \emptyset$, $\mb{u}^\ast \leftarrow \emptyset$, and $k \leftarrow 0$.
\State Solve LP relaxation of \eqref{eq:MILPAlg} with ADMM, Algorithm \ref{Alg:ADMM}. 
 Let $z^\text{LP}$ be its optimal value, and $\mb{\mu}^0$ be the dual values for the coupling constraints \eqref{eq:MILPAlg-Link}.
\State $lb\leftarrow z^\text{LP}$
\While{ some termination criteria is not met}
\State $k \leftarrow k+1$
\For{$\nu:=1$ to $N$}
\State solve $ \min\limits_{\mb{x}_\nu}  \{  {\cal L}_\nu(\mb{x}_\nu,\mb{\mu}^{k-1}): \mb{x}_\nu \in X_\nu\}$
\State let $v_\nu^k$ be the optimal value and $\mb{x}_\nu^{k}= (\mb{u}_\nu^{k},\mb{y}_\nu^{k})$ be an optimal solution
\EndFor
\If{$lb<  \mb{\mu}^\top \mb{b}+\sum\limits_{\nu\in {\cal P}}  v_\nu^k$}
\State $lb\leftarrow  \mb{\mu}^\top \mb{b}+\sum\limits_{\nu\in {\cal P}} v_\nu^k$
\EndIf
\State  $\mb{\mu}^{k} \leftarrow \mb{\mu}^{k-1}+\rho_\mu^{k}\left(\mb{b}-
\sum\limits_{\nu\in {\cal P} } \mb{A}_\nu  \mb{x}_\nu^k  \right) $ 
\If{$\mb{u}_\nu^{k} \notin {\cal S}$}
\State ${\cal S} \leftarrow {\cal S}\cup \{ \mb{u}_\nu^{k} \}$
\State compute $z(\mb{u}_\nu^{k})$ by solving \eqref{eq:LP1} with ADMM, Algorithm \ref{Alg:ADMM}
\If{$z(\mb{u}_\nu^{k+1})<ub$ }
\State $ub \leftarrow z(\mb{u}_\nu^{k})$
\State $\mb{u}^\ast \leftarrow \mb{u}_\nu^{k}$
\EndIf
\EndIf
\EndWhile
\end{algorithmic}
\end{algorithm}

\subsection{Release-and-Fix Heuristic}

\begin{algorithm}
\caption{Release-and-Fix Heuristic for MILPs \citep{Feizollahi:2015Large} }\label{Alg:RaF}
\begin{algorithmic}[1]
\State $ub \leftarrow +\infty$, $\mb{u}^\ast \leftarrow \emptyset$, and $k \leftarrow 0$.
\State {\bf ADMM-CR}: Solve LP relaxation of \eqref{eq:MILPAlg} with ADMM, Algorithm \ref{Alg:ADMM}. 
 Let $z^\text{LP}$ be its optimal value, and $\mb{\mu}^0$ be the dual values for the coupling constraints \eqref{eq:MILPAlg-Link}.
\State $lb\leftarrow z^\text{LP}$
\While{ time or iteration limits are not met}
\State $k \leftarrow k+1$
\State {\bf ADMM-Bin+}: Continue ADMM, Algorithm \ref{Alg:ADMM}, for the original MILP \eqref{eq:MILPAlg}  until some criteria are not met. In this phase, binary variables are restricted to take only 0 or 1 values. Let $\hat{\mb{u}}$ be the binary subvector of the current solution at the end of this phase.
\State {\bf ADMM-Bin-}: Fix the binary variables at their level of $\hat{\mb{u}}$. Continue ADMM, Algorithm \ref{Alg:ADMM}, to compute $z(\mb{\hat u})$ by solving LP \eqref{eq:LP1} 
\If{$z(\mb{\hat u})<ub$ }
\State $ub \leftarrow z(\mb{\hat u})$
\State $\mb{u}^\ast \leftarrow \mb{\hat u}$
\EndIf
\EndWhile
\end{algorithmic}
\end{algorithm}

\cite{Feizollahi:2015Large} have developed an ADMM based a heuristic decomposition method, which was called release-and-fix to solve MILPs arising from electric power network unit commitment problems. Algorithm \ref{Alg:RaF} presents a high level scheme of the release-and-fix method. This algorithm along with some refinements were able to mitigate oscillations and traps in local optimality. This method was able to find very good solutions with relatively small optimality gap for large scale unit commitment problems \citep{Feizollahi:2015Large}. But, it was not able to get the exact solution of MILP \eqref{eq:MILPAlg}.

\section{Exact Distributed Algorithms} \label{Sec:ExactAlgs}
In this section, we propose a distributed MILP approach where each block solves its own modified LR subproblem iteratively. 
The approach evaluates the cost of binary solutions as candidate partial solutions and refines them to get a primal feasible solutions to the overall problem. To improve the lower bound and prevent cycling in Algorithm \ref{Alg:DualD}, the explored binary solutions are then cut-off from future consideration in all subproblems. 

This idea is similar to the scenario decomposition algorithm for two-stage 0-1 stochastic MILP problems proposed in \citep{Ahmed:2013scenario}. In the two-stage 0-1 stochastic MILP model at  \citep{Ahmed:2013scenario}, each scenario is assumed to be a block and nonanticipativity constraints are coupling different scenarios.  In that model, binary variables are only present in the first stage and they are the same for different scenarios. Therefore, it is straightforward to cutoff  explored binary solutions from the feasible regions of all subproblems. On the contrary, in our loosely coupled MILP model \eqref{eq:MILPAlg}, binary variables are not the same for different blocks. Then, it is not clear how to cutoff a global binary solution from the feasible regions of subproblems. For instance, in Example \ref{Ex:MILPAlgExample}, $\mb{u}_1=(u_{11}, u_{12}, u_{13})^\top$ and  $\mb{u}_2=(u_{21}, u_{22}, u_{23})^\top$ are completely different binary vectors for blocks 1 and 2, respectively.  In Example \ref{Ex:MILPAlgExample}, consider $\hat{\mb u}=(\hat{\mb u}_1^\top,\hat{\mb u}_2^\top)^\top$ where $\hat{\mb u}_1=(1,1,0)^\top\in U_1$ and $\hat{\mb u}_2=(0,0,0)^\top\in U_2$. Then, it is a challenge to cutoff  $\hat{\mb u}=(\hat{\mb u}_1^\top,\hat{\mb u}_2^\top)^\top=(1,1,0,0,0,0)^\top \in U_1\times U_2$ from the local feasible regions of blocks 1 and 2 in a distributed and parallel fashion.
 In this section, we propose two exact algorithms to handle this process in a distributed framework. 
 
 For given $\mb{\mu}\in \Re^m$ and ${\cal S}\subset U$, we define the \emph{restricted Lagrangian relaxation} (RLR)
\begin{equation}\label{eq:MILP1-RLR}
\begin{split}
z^\text{RLR}(\mb{\mu}, {\cal S}):=\mb{\mu}^\top \mb{b}+ \min\limits_{ \mb{x}_1, \cdots , \mb{x}_N }  & \sum_{ \nu\in {\cal P} } {\cal L}_\nu(\mb{x}_\nu,\mb{\mu})
\\
\text{s.t. } &\mb{x}_\nu \in X_\nu,\, \forall \nu\in {\cal P}, \\
& \mb{u}\notin {\cal S}.
\end{split}
\end{equation}
Recall from Assumption \ref{Assump:MILPAlg}, $\mb{x}$ consists of the binary variables' subvector $\mb{u}$ and the continuous variables' subvector $\mb{y}$. Note that, $ub({\cal S}) :=\min\limits_{\mb{\hat u} \in {\cal S}} \{ z(\mb{\hat u}) \}$
and $lb(\mb{\mu}, {\cal S}):=\min \{ z^\text{RLR}(\mb{\mu}, {\cal S}), ub( {\cal S}) \}$ are valid upper and lower bounds for $z^\text{IP}$, respectively. 

\begin{proposition}\label{prop:finite-iter1}
Consider MILP \eqref{eq:MILPAlg} under Assumption \ref{Assump:MILPAlg}. For any $\mb{\mu}\in \Re^m$, there exists a set ${\cal S}\subset U$ such that $ub({\cal S}) =lb(\mb{\mu}, {\cal S})=z^\text{IP}$.
\end{proposition}

\proof{Proof}
By \eqref{eq:z-U}, we know that $ub(U) :=\min\limits_{\mb{\hat u} \in U} \{ z(\mb{\hat u}) \}=z^\text{IP}$.
Clearly, $z^\text{RLR}(\mb{\mu}, U)=+\infty$ and consequently $lb(\mb{\mu}, U):=\min \{ z^\text{RLR}(\mb{\mu}, U), ub(U) \}=ub(U)$. $\Box$
\endproof
 
Note that for any $\mb{\mu}\in \Re^m$, $z^\text{RLR}(\mb{\mu}, {\cal S})$ and $ub({\cal S})  $ are non-decreasing and non-increasing, respectively, functions of ${\cal S}$, 
 i.e.
$z^\text{RLR}(\mb{\mu}, {\cal S}) \le z^\text{RLR}(\mb{\mu}, {\cal T})$  and $ub( {\cal S}) \ge ub({\cal T})$ for any pair of sets ${\cal S}$ and ${\cal T}$ such that ${\cal S} \subset {\cal T} \subset U$. 
Therefore, for any $\mb{\mu}\in \Re^m$, there exists a set ${\cal S}(\mb{\mu} )\subset U$ such that $z^\text{RLR}(\mb{\mu}, {\cal S}(\mb{\mu} ) )\ge ub({\cal S}(\mb{\mu} ))$ and consequently $lb(\mb{\mu}, {\cal S}(\mb{\mu} ))=z^\text{IP}=ub({\cal S}(\mb{\mu} ))$. In other words, it is possible to close the duality gap for MILP \eqref{eq:MILPAlg} by cutting off some finite number of binary solutions in \eqref{eq:MILP1-RLR} via constraints $ \mb{u}\notin {\cal S}$.

For a given binary vector $\hat{\mb{u}}\in \{0,1\}^{n^1}$ let us define the simple binary cut (SBC) of $\hat{\mb{u}}$ in terms of binary decision vector $\mb{u}\in \{0,1\}^{n^1}$ as follows:
\begin{equation}\label{eq:SBC}
\text{SBC} (\mb{u},\mb{\hat u}):~~ \sum_{k:\hat{u}_k=0} u_k + \sum_{k:\hat{u}_k=1} (1-u_k) \ge 1.
\end{equation}
Then, $\text{SBC} (\mb{u},\mb{\hat u})$ for $\hat{\mb u}=(1,1,0,0,0,0)^\top $ in Example \ref{Ex:MILPAlgExample} is the following inequality: 
\begin{equation}\label{eq:SBC-Ex}
-u_{11}-u_{12}+u_{13}+u_{21}+u_{22}+u_{23} \ge -1.
\end{equation}
To cutoff multiple solutions, stronger cuts can be used as described in \citep{Angulo:2015}.
Using the concept of SBC, the constraint $\mb{u}\notin {\cal S}$ in \eqref{eq:MILP1-RLR} can be represented as $\text{SBC} (\mb{u},\mb{\hat u})$, for all $\mb{\hat u} \in {\cal S}$. However this constraint couples different blocks to each other  and defeats the goal of problem decomposition.
For example, in constraint \eqref{eq:SBC-Ex}, all binary variables from blocks 1 and 2 are present.
 Next, we propose different techniques to overcome this issue by introducing equivalent formulations of \eqref{eq:MILP1-RLR} which are decomposable.

\subsection{Binary Variables Duplication}
In our first approach of decoupling the constraint $\mb{u}\notin {\cal S}$ in \eqref{eq:MILP1-RLR}, we propose to duplicate the whole vector of binary variables and give a copy of it to each block.
For each pair of $\nu,\nu'\in {\cal P}$, let $\tilde{\mb u}_{\nu,\nu'}\in \tilde{U}_{\nu,\nu'} \subset \{0,1\}^{n_{\nu'}^1}$ be block $\nu$'s perception of $\mb{u}_{\nu'}$, where $U_{\nu,\nu'}$ is the set of all possible values for $\tilde{\mb u}_{\nu,\nu'}$. 
For convenience, let $\tilde{\mb u}_\nu$ and $\tilde{U}_{\nu}$ be block $\nu$'s perception of ${\mb u}$ and $U$.
Note that $\tilde{\mb u}_\nu \in \{0,1\}^{n^1}$ and $\tilde{U}_{\nu}\subset \{0,1\}^{n^1}$, for all $\nu\in {\cal P}$. 

It can be assumed $U_{\nu'} \subset \tilde{U}_{\nu,\nu'} $ for all $\nu\ne \nu'$ where it is possible that $U_{\nu'} \ne \tilde{U}_{\nu,\nu'}$. 
For example one may assume $\tilde{U}_{\nu,\nu'}= \{0,1\}^{n_{\nu'}^1}$. Therefore, it may happen $\tilde{U}_{\nu,\nu'} \backslash U_{\nu'} \ne \emptyset$; i.e. block $\nu$ may not know any explicit or implicit descriptions of $U_{\nu'}$ and consequently its perception of $\mb{u}_{\nu}$ can be infeasible. But, block $\nu$ should receive an infeasibility alert from block $\nu'$, if $\mb{\hat u}_{\nu,\nu'}\notin U_{\nu'}$.
Then, $\mb{\hat u}_{\nu,\nu'}$ can be cut off from $\tilde{U}_{\nu,\nu'}$ using SBC($\mb{u}_{\nu,\nu'},\mb{\hat u}_{\nu,\nu'}$) as defined in \eqref{eq:SBC}. 
In this algorithm, we assume $\tilde{U}_\nu=U$, for the sake of simplicity. Later, we will present other algorithms where the blocks do not need to know anything about the feasibility regions of the other blocks.

For Example \ref{Ex:MILPAlgExample},  blocks 1 and 2 perceptions of the overall binary vector $\mb{u}$ are
$\tilde{\mb u}_1=(\tilde{u}_{111},\tilde{u}_{112},\tilde{u}_{113},\tilde{u}_{121},\tilde{u}_{122},\tilde{u}_{123})^\top$ and 
$\tilde{\mb u}_2=(\tilde{u}_{211},\tilde{u}_{212},\tilde{u}_{213},\tilde{u}_{221},\tilde{u}_{222},\tilde{u}_{223})^\top$, respectively. In this case,
$\tilde{\mb u}_{11}=(\tilde{u}_{111},\tilde{u}_{112},\tilde{u}_{113})^\top \in U_1$, 
$\tilde{\mb u}_{12}=(\tilde{u}_{121},\tilde{u}_{122},\tilde{u}_{123})^\top\in U_1$, 
$\tilde{\mb u}_{21}=(\tilde{u}_{211},\tilde{u}_{212},\tilde{u}_{213})^\top\in U_1$, and
$\tilde{\mb u}_{22}=(\tilde{u}_{221},\tilde{u}_{222},\tilde{u}_{223})^\top\in U_2$. Then, $\text{SBC} (\mb{u},\mb{\hat u})$ cut \eqref{eq:SBC-Ex} for $\hat{\mb u}=(1,1,0,0,0,0)^\top $ can be reformulated as
\begin{equation}\label{eq:SBC-Ex1}
-\tilde{\mb u}_{111}-\tilde{\mb u}_{112}+\tilde{\mb u}_{113}+\tilde{\mb u}_{121}+\tilde{\mb u}_{122}+\tilde{\mb u}_{123} \ge 1,
\end{equation}
and 
\begin{equation}\label{eq:SBC-Ex2}
-\tilde{\mb u}_{211}-\tilde{\mb u}_{212}+\tilde{\mb u}_{213}+\tilde{\mb u}_{221}+\tilde{\mb u}_{222}+\tilde{\mb u}_{223} \ge 1.
\end{equation}
for blocks 1 and 2, respectively. Note that in inequality \eqref{eq:SBC-Ex1}, only (perception) binary variables from block 1 are present. Similarly, in inequality \eqref{eq:SBC-Ex2}, only (perception) binary variables from block 2 are present.

An equivalent formulation for \eqref{eq:MILP1-RLR} can be constructed by using the binary vectors $\tilde{\mb u}_1, \cdots, \tilde{\mb u}_N$, where all the blocks have the same perceptions of $\mb{u}$, i.e. 
\begin{equation}\label{eq:MILP3-Consensus}
\tilde{\mb u}_1= \cdots= \tilde{\mb u}_N,
\end{equation}
and the $\mb{u}\notin {\cal S}$ is replaced by
\begin{equation}\label{eq:MILP3-XS}
\tilde{\mb u}_\nu \in U \backslash {\cal S}.
\end{equation}
In Example \ref{Ex:MILPAlgExample}, constraint \eqref{eq:MILP3-Consensus} has the following form
\begin{equation}\nonumber
\begin{split}
\tilde{u}_{111}&=\tilde{u}_{211},\\
\tilde{u}_{112}&=\tilde{u}_{212},\\
\tilde{u}_{113}&=\tilde{u}_{213},\\
\tilde{u}_{121}&=\tilde{u}_{221},\\
\tilde{u}_{122}&=\tilde{u}_{222},\\
\tilde{u}_{123}&=\tilde{u}_{223}.
\end{split}
\end{equation}

Note that for all $\nu'\ne \nu$,  binary vectors $\tilde{\mb u}_{\nu, \nu'}$ are redundant. But, they make it possible to cut a global binary solution $\hat{\mb{u}}$ from the feasible region of all blocks. In other words, we use $\tilde{\mb u}_{\nu, \nu'}$ for $\nu' \ne \nu$  to handle constraint \eqref{eq:MILP3-XS}. 
Let $\tilde{\mb x}_\nu:=(\tilde{\mb u}_\nu,\mb{y}_\nu)\in \{0,1\}^{n^1}\times \Re^{n^2_\nu}$. Note that for all $\nu\in {\cal P}$, ${\mb u}_\nu$ is a subvector of $\tilde{\mb u}_\nu$ and consequently, $\mb{x}_\nu=({\mb u}_\nu,\mb{y}_\nu)$ is  a subvector of $\tilde{\mb x}_\nu$.
Then, problem \eqref{eq:MILP1-RLR} can be reformulated as follows:
\begin{equation}\label{eq:MILP1-RLR2}
\begin{split}
z^\text{RLR}(\mb{\mu}, {\cal S})=\mb{\mu}^\top \mb{b}+ \min\limits_{ \tilde{\mb x}_1,\cdots,\tilde{\mb x}_N }  & \sum_{ \nu\in {\cal P} } {\cal L}_\nu(\mb{x}_\nu,\mb{\mu})
\\
\text{s.t. } & \mb{x}_\nu \in X_\nu \text{ and } \tilde{\mb u}_\nu \in U \backslash {\cal S}, \, \forall \nu\in {\cal P}\\
& \tilde{\mb u}_1= \cdots= \tilde{\mb u}_N.
\end{split}
\end{equation}

In the model \eqref{eq:MILP1-RLR2}, the  consensus constraints \eqref{eq:MILP3-Consensus} are joint between different blocks. 
To decouple these constraints, we use vectors of dual variables $\mb{\lambda}_\nu \in \Re^{n^1}$, for all $\nu$ such that $\sum_{\nu\in{\cal P} } \mb{\lambda}_\nu= \mb{0}$. Then, the new restricted Lagrangian relaxation for the model \eqref{eq:MILPAlg} is
\begin{equation}\label{eq:MILP3-LR}
\begin{split}
z^{\text{RLR}'}(\mb{\mu},\mb{\lambda},{\cal S}):=\mb{\mu}^\top \mb{b}+ \min\limits_{ \tilde{\mb x}_1,\cdots,\tilde{\mb x}_N }  & \sum_{\nu\in {\cal P}} {\cal L}'_\nu(\tilde{\mb x}_\nu,\mb{\mu},\mb{\lambda}_\nu)\\ 
\text{s.t. } &  \mb{x}_\nu \in X_\nu \text{ and } \tilde{\mb u}_\nu \in U \backslash {\cal S}, \, \forall \nu \in {\cal P},
\end{split}
\end{equation}
where $\mb{\lambda}=(\mb{\lambda}_1^\top,\cdots, \mb{\lambda}_N^\top)^\top$ and 
$
{\cal L}'_\nu(\tilde{\mb x}_\nu,\mb{\mu},\mb{\lambda}_\nu):= (\mb{c}^\top_\nu-\mb{\mu}^\top \mb{A}_\nu) \mb{x}_\nu   +\mb{\lambda}_\nu^\top \tilde{\mb u}_\nu.
$
To solve problem \eqref{eq:MILP3-LR}, it is sufficient for each block $\nu$ to solve its subproblem of 
$\min\limits_{ \tilde{\mb x}_\nu } \{ {\cal L}'_\nu(\tilde{\mb x}_\nu,\mb{\mu},\mb{\lambda}_\nu): \mb{x}_\nu \in X_\nu \text{ and } \tilde{\mb u}_\nu \in U \backslash {\cal S}\}$.
Note that $z^{\text{RLR}'}(\mb{\mu},\mb{\lambda},{\cal S}) \le z^{\text{RLR}}(\mb{\mu},{\cal S})$, for all ${\cal S}\subset U$, $\mb{\mu}\in \Re^m$ and $\mb{\lambda}_\nu\in \Re^{n^1}$, $\forall \nu\in {\cal P}$ such that $\sum\limits_{\nu\in{\cal P} } \mb{\lambda}_\nu= \mb{0}$. Moreover, $z^{\text{RLR}'}(\mb{\mu},\mb{\lambda},{\cal S})$ is a non-decreasing function of ${\cal S}$.

\begin{algorithm}
\caption{Distributed MILP with Binary Variables Duplication}\label{Alg:DMILP1}
\begin{algorithmic}[1]
\State Run Algorithm \ref{Alg:DualD} to initialize $ub$, $lb$, $\mb{u}^\ast$, $\mb{\mu}^0$ and ${\cal S}$.
\State $\mb{\lambda}^0 \leftarrow \mb{0}$ and $k \leftarrow 0$.
\While{$ub>lb$}
\State Lower bounding:
\While{ some termination criteria is not met}\label{Alg1:LB-start}
\State $k \leftarrow k+1$
\For{$\nu:=1$ to $N$}
\State solve $ \min\limits_{ \tilde{\mb x}_\nu } \{ {\cal L}'_\nu(\tilde{\mb x}_\nu,\mb{\mu}^{k-1},\mb{\lambda}_\nu^{k-1}): \mb{x}_\nu \in X_\nu \text{ and } \tilde{\mb u}_\nu \in U \backslash {\cal S}\}$.
\State let $v_\nu^k$ be the optimal value and $\tilde{\mb x}_\nu^{k}=(\tilde{\mb u}_\nu^{k},\mb{y}_\nu^{k})$ be an optimal solution
\EndFor
\If{$lb<  \mb{\mu}^\top \mb{b}+\sum\limits_{ \nu \in {\cal P}} v_\nu^k$}
\State $lb\leftarrow  \min\left\{ub, \mb{\mu}^\top \mb{b}+\sum\limits_{ \nu \in {\cal P}} v_\nu^k \right\}$
\EndIf
\State $\bar{\mb u}^{k}\leftarrow \frac{1}{|{\cal P}|} \sum\limits_{\nu\in {\cal P} } \tilde{\mb u}^{k}_{\nu}$
\State $\mb{\mu}^{k} \leftarrow \mb{\mu}^{k-1}+\rho^{k}_\mu \left(\mb{b}-
\sum\limits_{\nu\in {\cal P} } \mb{A}_\nu  \mb{x}_\nu^k  \right) $
 and  $\mb{\lambda}^{k}_\nu \leftarrow \mb{\lambda}^{k-1}_\nu+\rho^{k}_\lambda \left(\tilde{\mb u}^{k}_\nu-\bar{\mb u}^{k}  \right) $ 
\EndWhile \label{Alg1:LB-finish}
\State Let  $\hat{{\cal S}}^k=\cup_{\nu \in {\cal P}} \{\tilde{\mb u}_\nu^{k}\}$.
\State Upper bounding: run Algorithm \ref{Alg:UpperBounding} for set $\hat{{\cal S}}$ to update $ub$ and $\mb{u}^\ast$.
\State $ {\cal S} \leftarrow  {\cal S}  \cup \hat{{\cal S}}^k $
\EndWhile
\end{algorithmic}
\end{algorithm}

Let $\rho_\mu^k, \rho_\lambda^k>0$ be the step size for updating the dual vectors $\mb{\mu}$ and $\mb{\lambda}$ at iteration $k$. 
Then, our first exact distributed MILP method is as Algorithm \ref{Alg:DMILP1}. This algorithm is initialized by running ADMM to solve the LP relaxation and then switches to dual decomposition. In fact, this step initializes upper and lower bounds as well as dual vectors.
In the lower bounding loop (lines \ref{Alg1:LB-start}-\ref{Alg1:LB-finish}) of Algorithm \ref{Alg:DMILP1}, problem \eqref{eq:MILP3-LR} is solved in parallel by each block and the dual vectors $\mb{\mu}$ and $\mb{\lambda}$ are updated as well as the lower bound and candidate binary subvectors. Then, each candidate binary subvector is evaluated by solving an LP with ADMM method. In this step, the upper bound is updated. Finally, the candidate binary subvectors are added to set ${\cal S}$ and consequently are cutoff from feasible regions of all blocks. The algorithm continues until the lower bound hits the upper bound.

\begin{proposition}\label{prop:finite-DMILP1}
Algorithm \ref{Alg:DMILP1} can find an optimal solution of  MILP \eqref{eq:MILPAlg} under Assumption \ref{Assump:MILPAlg} in a finite number of iterations. 
\end{proposition}

\proof{Proof}
In the worst case, Algorithm \ref{Alg:DMILP1} needs to be run until cutting off all binary solutions in $U$, which are finite. But  for any feasible dual vectors $\mb{\mu}$ and $\mb{\lambda}$, we know that $z^{\text{RLR}'}(\mb{\mu},\mb{\lambda},U)=+\infty> z^\text{IP}$ which implies $ub \ngtr lb$ and the algorithm terminates. $\Box$
\endproof

\subsection{Auxiliary Binary Variables}
In Algorithm \ref{Alg:DMILP1}, each block has as many binary variables as $n^1$, the number of overall binaries in the original MILP problem \eqref{eq:MILPAlg}. Moreover, each block $\nu$ needs to know the constraints defining the set $U_{\nu'}$,  for all $\nu'\ne \nu$ or to be able to check the feasibility of $\tilde{\mb u}_{\nu,\nu'}$. Next, we propose another algorithm by introducing some auxiliary binary variables, in which different blocks do not need to know about other blocks' binary variables or feasible regions.

For a given ${\cal S}\subset U$, let ${\cal S}_\nu$, for all $ \nu\in {\cal P}$, be the minimal sets such that ${\cal S}_\nu \subset U_\nu$  and ${\cal S}\subset {\cal S}_1 \times \cdots \times {\cal S}_N$. That is for all $\mb{\hat u}_\nu\in {\cal S}_\nu$ and $\nu \in {\cal P}$, there exists a $\mb{\hat u}\in {\cal S}$ such that the $\nu$th block of  $\mb{\hat u}$ is $\mb{\hat u}_\nu$. 
Let $S_\nu:=\{1,\cdots, |{\mathcal S}_\nu|\}$ and denote the $l$th solution of ${\cal S}_\nu$ by $\hat{\mb u}_\nu (l)$.

\begin{example} \label{Ex:MILPAlgExample2}
Consider Example \ref{Ex:MILPAlgExample} with ${\cal S}=\{(1,1,0,0,0,0),(1,1,0,0,1,1)\}$. Then, it holds ${\cal S}_1=\{(1,1,0)\}$ and ${\cal S}_2=\{(0,0,0),(0,1,1)\}$.
\end{example}

For $\nu,\nu'\in {\cal  P}$ and $l\in S_{\nu'}$, let $w_{\nu,\nu',l}$ be a binary variable which is $1$, if block $\nu$'s perception of $\mb{u}_{\nu'}$ is $\hat{\mb u}_{\nu'} (l)$, and $0$ otherwise. 
For convenience, let $w_{\nu,\nu',0}$ be a binary variable which is $1$, if block $\nu$'s perception of $\mb{u}_{\nu'}$ is not in ${\cal S}_{\nu'}$, and $0$ otherwise.  Then,
\begin{equation}\label{eq:binaryu}
w_{\nu,\nu',l}\in \{0,1\}, \, \forall \nu'\in N,\, l\in S_{\nu'}\cup \{0\}.
\end{equation}

Then, for Example \ref{Ex:MILPAlgExample2}, block 1 has auxiliary binary variables $w_{1,1,0}$, $w_{1,1,1}$, $w_{1,2,0}$, $ w_{1,2,1}$,  $w_{1,2,2}$.
Binary variable $w_{111}$ is 1  if and only if block 1 perception of $\mb{u}_1$ are $(1,1,0)$. Binary variables $w_{121}$ and $w_{122}$ are 1 if and only if blocks 1 perceptions of $\mb{u}_2$ are $(0,0,0)$ and $(0,1,1)$, respectively. Similarly, $w_{110}$ and $w_{120}$ are 1 if and only if blocks 1 perceptions of $\mb{u}_1$ and $\mb{u}_2$ do not exist in ${\cal S}_1$ and ${\cal S}_1$, respectively. Likewise, block 2 has auxiliary binary variables $w_{2,1,0}$, $w_{2,1,1}$, $w_{2,2,0}$,   $w_{2,2,1}$, $w_{2,2,2}$.

Note that block $\nu$ does not know the length of $\mb{u}_{\nu'}$ or the values in the $\hat{\mb u}_\nu (l)$, unless $\nu=\nu'$. Therefore,  $\mb{u}_\nu= \hat{\mb u}_\nu (l)$ if and only if $w_{\nu,\nu,l}=1$. This relation between the binary vector $\mb{u}_\nu$ and the binary variable $w_{\nu,\nu,l}$ can be imposed by constraints \eqref{eq:XeqXhat} and \eqref{eq:XnoteqXhat}. 
\begin{equation}\label{eq:XeqXhat}
 \left\{ \begin{array}{l l}
u_{\nu k} \ge w_{\nu,\nu,l}, &  \text{ if } \hat{u}_{\nu k}(l)=1 \\
u_{\nu k} \le 1-w_{\nu,\nu,l},  & \text{ Otherwise}
\end{array} \right.~~ \forall l\in S_{\nu}, k=1,\cdots,n_\nu^1,
\end{equation}
\begin{equation}\label{eq:XnoteqXhat}
\sum_{k:\hat{u}_{\nu k}(l)=0} u_{\nu k} + \sum_{k:\hat{u}_{\nu k}(l)=1} (1-u_{\nu k}) \ge w_{\nu, \nu,0}, \,  \forall   l\in S_\nu.
\end{equation}
Each block $\nu$ should consider exactly one of the binary solutions $\hat{\mb u}_{\nu'}$ in ${\cal S}_{\nu'}$, for all $\nu' \in {\cal P}$, i.e.
\begin{equation}\label{eq:ExactlyOne}
\sum_{l\in  S_{\nu}\cup \{0\} }w_{\nu,\nu',l}=1, \, \forall \nu'\in {\cal P}.
\end{equation}
Inequality \eqref{eq:GlobalCut0} cuts the explored binary solutions to prevent cycling. 
\begin{equation} \label{eq:GlobalCut0}
\sum_{\nu' \in {\cal P}} \left[\sum_{l:\mb{\hat u}_{\nu'}(l)\ne \mb{\hat u}_{\nu'}(s) } w_{\nu,\nu',l} + \sum_{l:\mb{\hat u}_{\nu'}(l)= \mb{\hat u}_{\nu'}(s)} (1-w_{\nu,\nu',l}) \right] \ge 1, \, \forall s\in {\mathcal S}, 
\end{equation}
Because of the constraints \eqref{eq:binaryu} and  \eqref{eq:ExactlyOne}, constraint \eqref{eq:GlobalCut0} can be strengthened as follows:
\begin{equation} \label{eq:GlobalCut}
\sum_{
\substack{
\nu'\in{\cal P}\\ l\in S_{\nu'}:\mb{\hat u}_{\nu'}(l)\ne \mb{\hat u}_{\nu'}(s) }}  w_{\nu,\nu',l}  \le N-1, \, \forall s\in {\mathcal S}.
\end{equation}

Constraints \eqref{eq:XeqXhat}-\eqref{eq:ExactlyOne}, and \eqref{eq:GlobalCut} for block 2 in Example \ref{Ex:MILPAlgExample2} have the following form:
\begin{equation}\nonumber
\begin{split}
& \left.\begin{array}{l}
u_{21}\le 1- w_{221}, ~ u_{22}\le 1- w_{221}, ~u_{23}\le 1- w_{221},\\
u_{21}\le 1- w_{222}, ~ u_{22}\ge w_{222}, ~ u_{23}\ge w_{222},
\end{array} \right\} \text{Constraint \eqref{eq:XeqXhat}}\\
\end{split}
\end{equation}
\begin{equation}\nonumber
\begin{split}
& \left.\begin{array}{l}
 u_{21}+u_{22}+u_{22} \ge w_{220},\\
u_{21}+1-u_{22}+1-u_{22} \ge w_{220},\\
\end{array} \right\} \text{Constraint \eqref{eq:XnoteqXhat}}\\
\end{split}
\end{equation}
\begin{equation}\nonumber
\begin{split}
& \left.\begin{array}{l}
w_{210}+w_{211}=1,\\
 w_{220}+w_{221}+w_{223}=1,
\end{array} \right\} \text{Constraint \eqref{eq:ExactlyOne}}\\
\end{split}
\end{equation}
\begin{equation}\nonumber
\begin{split}
& \left.\begin{array}{l}
 w_{211}+w_{221}\le 1,\\
 w_{211}+w_{222}\le 1.
\end{array} \right\} \text{Constraint \eqref{eq:GlobalCut}}\\
\end{split}
\end{equation}

Let $\mb{w}_\nu$ be the vector of all binary variables $w_{\nu,\nu',l}$, for all $\nu'\in {\cal P}$ and all $l\in S_{\nu'}$
In the second distributed MILP algorithm, we use the auxiliary binary vector $\mb{w}_\nu\in \{0,1\}^{|{\cal P}|+\sum\limits_{\nu' \in {\cal P}} |{\cal S}_{\nu'}| }$, for all $\nu \in {\cal P}$, to develop another equivalent model for \eqref{eq:MILP1-RLR}. Considering the  consensus constraints
\begin{equation}\label{eq:ConsensusU}
 \mb{w}_1=\cdots= \mb{w}_N,
\end{equation}
problem \eqref{eq:MILP1-RLR} can be reformulated as follows.
\begin{equation}\label{eq:MILP-Alg2}
\begin{split}
z^\text{RLR}(\mb{\mu},{\cal S})=\mb{\mu}^\top \mb{b}+ \min\limits_{\mb{x},  \mb{w}_1,\cdots,\mb{w}_N }  & \sum_{ \nu\in {\cal P} } {\cal L}_\nu(\mb{x}_\nu,\mb{\mu})\\
\text{s.t. } &   \mb{x}_\nu \in X_\nu \text{ and } \eqref{eq:binaryu}-\eqref{eq:ExactlyOne},\eqref{eq:GlobalCut}, \,  \forall \nu\in {\cal P}, \\
&   \mb{w}_1=\cdots= \mb{w}_N.
\end{split}
\end{equation}
Consensus constraints \eqref{eq:ConsensusU} are coupling different block in the problem \eqref{eq:MILP-Alg2}.
To decouple these constraints, we use the feasible dual variable vectors $\mb{\gamma}_\nu \in \Re^{{|{\cal P}|+\sum\limits_{\nu' \in {\cal P}} |{\cal S}_{\nu'}| }}$, for all $ \nu \in {\cal P}$ such that $\sum_{\nu\in{\cal P} } \mb{\gamma}_\nu= \mb{0}$. Then, the new restricted Lagrangian relaxation for the model \eqref{eq:MILPAlg} is
\begin{equation}\label{eq:MILP-RLR3}
\begin{split}
z^{\text{RLR}''} (\mb{\mu},\mb{\gamma},{\cal S}):=\mb{\mu}^\top \mb{b}+ & \min\limits_{\mb{x},  \mb{w}_1,\cdots,\mb{w}_N }   \sum_{\nu\in {\cal P}} {\cal L}''_\nu(\mb{x}_\nu,\mb{w}_\nu,\mb{\mu},\mb{\gamma}_\nu)\\ 
&~~ \text{s.t. } ~~  \mb{x}_\nu \in X_\nu, \text{ and } \eqref{eq:binaryu}-\eqref{eq:ExactlyOne},\eqref{eq:GlobalCut}, \,  \forall i\in {\cal N},
\end{split}
\end{equation}
where $\mb{\gamma}=(\mb{\gamma}_1, \cdots, \mb{\gamma}_N)$ and ${\cal L}''_\nu(\mb{x}_\nu,\mb{w}_\nu,\mb{\mu},\mb{\gamma}_\nu):= (\mb{c}^\top_\nu-\mb{\mu}^\top \mb{A}_\nu) \mb{x}_\nu +\mb{\gamma}_\nu \mb{w}_\nu$. 
Note that $z^{\text{RLR}''} (\mb{\mu},\mb{\gamma},{\cal S}) \le z^{\text{RLR}}(\mb{\mu},{\cal S})$, for all ${\cal S}\subset U$, and feasible dual variable vectors  $\mb{\mu}$ and $\mb{\gamma}$. Moreover, $z^{\text{RLR}''} (\mb{\mu},\mb{\gamma},{\cal S})$ is a non-decreasing function of ${\cal S}$.

Let $\rho_\gamma^k>0$ be the step size for updating the dual vector $\mb{\gamma}$ at iteration $k$. Then, our second exact distributed MILP approach is as Algorithm \ref{Alg:DMILP2}. The overall scheme of Algorithm \ref{Alg:DMILP2} is similar to Algorithm \ref{Alg:DMILP1}. The main difference is that instead of problem \eqref{eq:MILP3-LR}, problem  \eqref{eq:MILP-RLR3} is solved in parallel in
the lower bounding loop (lines \ref{Alg2:LB-start}-\ref{Alg2:LB-finish}) of Algorithm \ref{Alg:DMILP2}. Different blocks do not need to know about other blocks' vector $\mb{u}_\nu$ of binary variables or feasible regions $U_\nu$ to solve problem  \eqref{eq:MILP-RLR3} in parallel. Moreover, in line \ref{Step:1} of Algorithm \ref{Alg:DMILP2}, a new binary solution is added to ${\cal S}_\nu$ which results in adding a new corresponding binary variable $w$ and a new dual variable $\gamma$ to all blocks.

\begin{algorithm}
\caption{Distributed MILP with Auxiliary Binary Variables}\label{Alg:DMILP2}
\begin{algorithmic}[1]
\State Run Algorithm \ref{Alg:DualD} to initialize $ub$, $lb$, $\mb{u}^\ast$, $\mb{\mu}^0$ and ${\cal S}$.
\State Based on ${\cal S}$, set up the sets ${\cal S}_\nu$, for all $\nu\in {\cal P}$.
\State $\mb{\gamma}^0 \leftarrow \mb{0}$ and $k \leftarrow 0$.
\While{$ub>lb$}
\State Lower bounding:
\While{ some termination criteria is not met} \label{Alg2:LB-start}
\State $k \leftarrow k+1$.
\For{$\nu:=1$ to $N$}
\State solve $ \min\limits_{\mb{x}_\nu,  \mb{w}_\nu } \{{\cal L}''_\nu(\mb{x}_\nu,\mb{w}_\nu,\mb{\mu}^{k-1},\mb{\gamma}_\nu^{k-1})
:  \mb{x}_\nu \in X_\nu,  \eqref{eq:binaryu}-\eqref{eq:ExactlyOne}, \eqref{eq:GlobalCut}\}$
\State let $v_\nu^k$ be the optimal value and $(\mb{x}_\nu^{k},  \mb{w}_\nu^{k})$ be an optimal solution
\EndFor
\If{$lb<  \mb{\mu}^\top \mb{b}+\sum\limits_{\nu \in {\cal P}} v_\nu^k$}
\State $lb\leftarrow \min\left\{ub, \mb{\mu}^\top \mb{b}+\sum\limits_{\nu \in {\cal P}} v_\nu^k \right\}$
\EndIf
\State $\bar{\mb w}^{k}\leftarrow \frac{1}{|{\cal P}|} \sum\limits_{\nu\in {\cal P} } {\mb w}^{k}_{\nu}$
\State $\mb{\mu}^{k} \leftarrow \mb{\mu}^{k-1}+\rho^{k}_\mu \left(\mb{b}-\sum\limits_{\nu\in {\cal P} } \mb{A}_\nu  \mb{x}_\nu^k  \right) $ and  $\mb{\gamma}^{k}_\nu \leftarrow \mb{\gamma}^{k-1}_\nu+\rho^{k}_\gamma \left({\mb w}^{k}_\nu-\bar{\mb w}^{k}  \right) $ 
\EndWhile\label{Alg2:LB-finish}
\For{$\nu:=1$ to $N$}
\If{$\sum_{\nu' \in  {\cal P}} w_{\nu,\nu',0} \ge 1$}
\State ${\cal S}_\nu \leftarrow {\cal S}_\nu  \cup  \{{\mb u}_\nu (0)\}$ \label{Step:1}
\EndIf
\EndFor
\State Let $\tilde{\mb u}_\nu^{k}$ be the corresponding $\tilde{\mb u}_\nu \in U$ to $\mb{w}_\nu^{k}$
\State   $\hat{{\cal S}}\leftarrow \cup_{\nu \in {\cal P}} \{\tilde{\mb u}_\nu^{k}\}$. 
\State Upper bounding: run Algorithm \ref{Alg:UpperBounding} for set $\hat{{\cal S}}$ to update $ub$ and $\mb{u}^\ast$.
\State $ {\cal S} \leftarrow  {\cal S}  \cup \hat{{\cal S}}^k $
\EndWhile
\end{algorithmic}
\end{algorithm}

\begin{proposition}\label{prop:finite-DMILP2}
Algorithm \ref{Alg:DMILP2} can find an optimal solution of  MILP \eqref{eq:MILPAlg} under Assumption \ref{Assump:MILPAlg} in a finite number of iterations. 
\end{proposition}

\proof{Proof}
In the worst case, Algorithm \ref{Alg:DMILP2} needs to be run until cutting off all binary solutions in $U$, which are finite. But  for any feasible dual vectors $\mb{\mu}$ and $\mb{\gamma}$, we know that $z^{\text{RLR}''}(\mb{\mu},\mb{\gamma},U)=+\infty> z^\text{IP}$ which implies $ub \ngtr lb$ and the algorithm terminates. $\Box$
\endproof

\section{Illustrative Computations}\label{Sec:Computations}

In this section, we present numerical results testing the exact distributed MILP Algorithms \ref{Alg:DMILP1} and \ref{Alg:DMILP2} presented in Section \ref{Sec:ExactAlgs}, on small UC instances. We used 6 small unit commitment (UC) instances with 3, 4 and 5 generators for $T$=12 and 24 hours of planning. For details of UC formulation which is a MILP problem the reader can see \citep{Carrion:2006, Feizollahi:2015Large, Costley:2017}.
Table \ref{table:ExactTestCases} presents details of these instances. In Table \ref{table:ExactTestCases}, ``\# Gen'' and  ``Gen. types'' denote the number and types of generator in each instance (see Table  \ref{Table:GeneratorData} for details of each generator type). The total system demand at each hour is determined as given in Table \ref{Table:TotalDemand}.
The labels ``\# Bin. Vars.'', ``\# Cont. Vars.'', and ``\# Constr.'' denote the number of binary variables, continuous variables, and constraints, respectively, for each test case.
Moreover, the columns
$z^\text{LP}$, $z^\text{IP}$, ``Duality Gap'', and $t_\text{C}$ represent optimal objective value of LP relaxation and MILP formulation for UC, relative duality gap in percentage (between $z^\text{LP}$ and $z^\text{IP}$), and the solution time (in seconds) in central approach, respectively. An estimation for Lagrangian dual, which is obtained as the best lower bound in 100 iterations of the dual decomposition method, is denoted by $\tilde{z}^\text{LD}$. Note that finding an optimal vector of dual variables in the dual decomposition algorithm is not guaranteed. Then,  $\tilde{z}^\text{LD}$ is not necessarily equal or close to the value of  Lagrangian dual.

\begin{table}
\caption{Generator Data \citep{Carrion:2006}}
\label{Table:GeneratorData}
\center
\setlength\extrarowheight{2pt}
\scalebox{1}{
\begin{tabular}{|@{\hspace{0.2mm}}c@{\hspace{0.2mm}} |@{\hspace{0.2mm}}c@{\hspace{0.2mm}} |@{\hspace{0.2mm}}c@{\hspace{0.2mm}} |@{\hspace{0.2mm}}c@{\hspace{0.2mm}} |@{\hspace{0.2mm}}c@{\hspace{0.2mm}} |@{\hspace{0.2mm}}c@{\hspace{0.2mm}} |@{\hspace{0.2mm}}c@{\hspace{0.2mm}} |@{\hspace{0.2mm}}c@{\hspace{0.2mm}} |@{\hspace{0.2mm}}c@{\hspace{0.2mm}} |@{\hspace{0.2mm}}c@{\hspace{0.2mm}} |@{\hspace{0.2mm}}c@{\hspace{0.2mm}} |@{\hspace{0.2mm}}c@{\hspace{0.2mm}} |}
\hline
\multirow{3}{*}{Gen} & \multicolumn{6}{c|}{Technical Information} & \multicolumn{5}{c|}{Cost Coefficients}\\ \cline{2-12} 
&  \multicolumn{1}{c|}{$\overline{P}$} & \multicolumn{1}{c|}{$\underline{P}$}	& \scriptsize{ TU/TD }& \scriptsize{RU/RD} & $T^{\text{Init}}$ & $T^{\text{cold}}$	& $C^{\text{NL}}$	&	$C^{\text{LV}}$ & $C^{\text{Q}}$ & $C^{\text{HS}}$  & $C^{\text{CS}}$ \\
& (MW)	&	(MW)	&	(h)	&	(MW/h)	&	(h)	&	(h)	&	(\$/h)	& (\$/MWh)	& (\$/MW$^2$h)	&	(\$)	&	(\$) \\ 
\hhline{|============|}
1 & 455 & 150 & 8 & 225 & +8 & 5 & 1000 & 16.19 & 0.00048 & 4500 & 9000 \\[1pt] \hline 
2 & 455 & 150 & 8 & 225 & +8 & 5 & 970 & 17.26 & 0.00031 & 5000 & 10000 \\[1pt] \hline 
3 & 130 & 20 & 5 & 50 & -5 & 4 & 700 & 16.60 & 0.00200 & 550 & 1100 \\[1pt] \hline 
4 & 130 & 20 & 5 & 50 & -5 & 4 & 680 & 16.50 & 0.00211 & 560 & 1120 \\[1pt] \hline 
5 & 162 & 25 & 6 & 60 & -5 & 4 & 450 & 19.70 & 0.00398 & 900 & 1800 \\[1pt] \hline 
6 & 80 & 20 & 3 & 60 & -3 & 2 & 370 & 22.26 & 0.00712 & 170 & 340 \\[1pt] \hline 
7 & 85 & 25 & 3 & 60 & -3 & 2 & 480 & 27.74 & 0.00079 & 260 & 520 \\[1pt] \hline 
8 & 55 & 10 & 1 & 135 & -1 & 0 & 660 & 25.92 & 0.00413 & 30 & 60 \\[1pt] \hline 
\end{tabular}
}
\end{table}

\begin{table}
\caption{Total Demand (\% of Total Capacity)}
\label{Table:TotalDemand}
\center
\scalebox{.9}{
\begin{tabular}{|c|c|c|c|c|c|c|c|c|c|c|c|c|}
\hline
Time & 1 & 2 & 3 & 4 & 5 & 6 & 7 & 8 & 9 & 10 & 11 & 12 \\  
Demand & 71\% & 65\% & 62\% & 60\% & 58\% & 58\% & 60\% & 64\% & 73\% & 80\% & 82\% & 83\% \\    \hline
Time & 13 & 14 & 15 & 16 & 17 & 18 & 19 & 20 & 21 & 22 & 23 & 24 \\  
Demand & 82\% & 80\% & 79\% & 79\% & 83\% & 91\% & 90\% & 88\% & 85\% & 84\% & 79\% & 74\% \\  
  \hline
\end{tabular}
}
\end{table}

\begin{table}
\center
\caption{Test case details for exact algorithms}\label{table:ExactTestCases}
\scalebox{0.85}{
\begin{tabular}{|c|p{1.5cm}|c|p{1.2cm}|p{1.4cm}|c|c|c|p{1.5cm}|p{.8cm}|c|}
\hline 
\# Gen & Gen. types & T & \# Bin. Vars. & \# Cont. Vars. & \# Constr. & $z^\text{LP}$ & $z^\text{IP}$ & Duality Gap (\%) & $t_\text{C}$ (Sec) & $\tilde{z}^\text{LD}$\\ \hline \hline
\multirow{2}{*}{3}& \multirow{2}{*}{6,7,8} & 24 & 216 & 144 & 891 & 139896 & 146403 & 4.44 & 0.03 & 139933\\ \cline{3-11} 
& &  12 & 108 & 72 & 435 & 68212 & 70945 & 3.85 & 0.09 & 68226\\ \hline \hline
\multirow{2}{*}{4}& \multirow{2}{*}{3,5,6,8} & 24 & 288 & 192 & 1244 & 207068 & 212771 & 2.68 & 0.14 & 207100\\ \cline{3-11}
 &  & 12 & 144 & 96 & 596 & 101676 & 104381 & 2.59 & 0.09 & 101686\\ \hline \hline
\multirow{2}{*}{5}& \multirow{2}{*}{1,5,6,7,8} & 24 & 360 & 240 & 1514 & 354684 & 359197 & 1.26 & 0.22 & 354705\\ \cline{3-11}
 &  & 12 & 180 & 120 & 722 & 171099 & 172994 & 1.10 & 0.11 & 171110\\ \hline
\end{tabular}
}
\end{table}

All algorithms were coded in C++ using CPLEX 12.6 through the Concert API. Central UC instances were solved using internal CPLEX multi-threading with four cores.
The step sizes $\rho_\mu$, $\rho_\lambda$ and $\rho_\gamma$ were set to be 0.01, 10 and 50, respectively. The algorithms start with running ADMM to solve the LP relaxations of the UC instances to initialize the vector of dual variables $\mb{\mu}$ and the lower bound $lb$. Then, they do 100 iterations of the dual decomposition algorithm to improve the lower bound. Then, the main body of Algorithms \ref{Alg:DMILP1} and \ref{Alg:DMILP2} starts with 200 iterations limit where the first 10 iterations are spent on updating dual vectors $\mb{\lambda}$ and $\mb{\gamma}$ without adding cuts. In each iteration, the lower bounding phase does 10 sub-iterations. Then, new candidate binary vectors are explored by the upper bounding procedure and cutoff from the feasible regions of all blocks.

\begin{table}
\center
\caption{Summary of the results for the exact Algorithm \ref{Alg:DMILP1}}\label{table:ResultsExactAlg1}
\begin{tabular}{|c|c|r|r
|r|r|r|r|r|r|r|}
\hline
\# Gen & T & $t_\text{0}$ & $t_\text{1}$ & $t^\ast$ & $t_\text{all}$ & iter$_1$ & iter$^\ast$ & iter$_\text{all}$ & \# Feas. & \# Cut\\ \hline \hline
\multirow{2}{*}{3} & 24 & 3.85 & 4.02 & 4.67 & 4.69 & 1 & 5 & 5 & 12 & 16\\ \cline{2-11}
& 12 & 2.16 & 2.23 & 2.59 & 3 & 1 & 4 & 7 & 12 & 19\\ \hline \hline
\multirow{2}{*}{4} & 24 & 4.54 & 4.7 & 7.4 & 193.1 & 1 & 11 & 118 & 90 & 530\\ \cline{2-11}
& 12 & 2.23 & 2.31 & 5.36 & 34.9 & 1 & 16 & 61 & 38 & 252\\ \hline \hline
\multirow{2}{*}{5} & 24 & 5.07 & 5.29 & 97.08 & 1621.72 & 1 & 42 & 190$^\ast$ & 303 & 1004\\ \cline{2-11}
& 12 & 2.18 & 2.29 & 3.01 & 715.88 & 1 & 5 & 190$^\ast$ & 290 & 962\\ \hline
\end{tabular}
\end{table}

\begin{table}
\center
\caption{Summary of the results for the exact Algorithm \ref{Alg:DMILP2}}\label{table:ResultsExactAlg2}
\begin{tabular}{|c|c|r|r
|r|r|r|r|r|r|r|}
\hline
\# Gen & T & $t_\text{0}$ & $t_\text{1}$ & $t^\ast$ & $t_\text{all}$ & iter$_1$ & iter$^\ast$ & iter$_\text{all}$ & \# Feas. & \# Cut\\ \hline \hline
\multirow{2}{*}{3} & 24 & 4.07 & 4.39 & 4.39 & 13.47 & 2 & 2 & 45 & 4 & 9\\ \cline{2-11}
& 12 & 2.35 & 2.49 & 2.49 & 3.97 & 1 & 1 & 11 & 4 & 7\\ \hline \hline

\multirow{2}{*}{4} & 24 & 4.56 & 4.81 & 4.81 & 6.81 & 2 & 2 & 10 & 12 & 23\\ \cline{2-11}
& 12 &  2.07 & 2.17 & 2.17 & 2.87 & 1 & 1 & 6 & 6 & 11\\ \hline \hline
\multirow{2}{*}{5} & 24 & 4.45 & 5.05 & 5.05 & 601.12 & 4 & 4 & 190$^\ast$ & 84 & 373\\ \cline{2-11}
& 12 &  1.98 & 2.19 & 2.19 & 24.06 & 3 & 3 & 38 & 54 & 122\\ \hline
\end{tabular}
\end{table}

Summary of the results for exact Algorithms \ref{Alg:DMILP1} and \ref{Alg:DMILP2} are presented in Tables \ref{table:ResultsExactAlg1} and \ref{table:ResultsExactAlg2}, respectively.
In Tables \ref{table:ResultsExactAlg1} and \ref{table:ResultsExactAlg2}, $t_\text{0}$, $t_\text{1}$, $t^\ast$, and $t_\text{all}$ are the estimated parallel times spent to initialize the algorithm, to find the first and best feasible solution, and to terminate the algorithm, respectively.
The exact algorithms were initialized by running ADMM for the LP relaxation and 100 iterations of the dual decomposition.
``iter$_1$'', ``iter$^\ast$'', ``iter$_\text{all}$'' are the corresponding number of iteration to $t_\text{1}$, $t^\ast$, and $t_\text{all}$, respectively.
``\# Feas.'', ``\# Cut'' are the number of feasible explored solutions and cuts (all explored binary solution), respectively.

For the 5 generator cases with $T$=24 and 12,  Algorithm \ref{Alg:DMILP1} terminated with \% 1.078 and \%0.911 optimality gaps after 190 iterations.  For the 5 generator case with $T$=24, Algorithm \ref{Alg:DMILP2} terminated with \%0.671 optimality gap after 190 iterations.
All other cases were solved to optimality. Based on the results in Tables \ref{table:ResultsExactAlg1} and \ref{table:ResultsExactAlg2}, for most cases, Algorithm \ref{Alg:DMILP2} outperforms Algorithm \ref{Alg:DMILP1}, in the sense that it requires less solution time ($t_\text{all}$), total number of iterations (iter$_\text{all}$) and cuts.

\section{Conclusions and Future Work}\label{Sec:Conclusion}

In this paper, we proposed exact distributed algorithms to solve MILP problems. A key challenge is that, because of
the non-convex nature of MILPs, classical distributed and decentralized optimization
approaches cannot be applied directly to find their optimal solutions. 
The main contributions of the paper are as follows:
\begin{enumerate}
\item two exact distributed MILP algorithms which are able to optimally solve MILP problems in a distributed manner and output primal feasible solutions
\item  primal cuts were added to restrict the Lagrangian relaxation 
and improve the lower bound on the objective function of the original MILP problem.
\item illustrative computation on unit commitment problem. 
\end{enumerate}

The main conclusions are as follows:
\begin{enumerate}
\item The proposed exact algorithms are proof-of-concept implementations to verify possibility of obtaining the global optimal solutions of MILPs in a distributed manner. Hence, the focus is not on computational times or number of iterations. 

\item Algorithm \ref{Alg:DMILP2} requires less information exchange between block than Algorithm  \ref{Alg:DMILP1}.

\item Based on the results in Tables \ref{table:ExactTestCases}-\ref{table:ResultsExactAlg2}, these exact distributed algorithms take much more time than the central approach. In particular, the solution times for Algorithms \ref{Alg:DMILP1} and \ref{Alg:DMILP2} are 3 seconds to 30 minutes  while the central problems are solved in less than a second. 

\item In general, Algorithm  \ref{Alg:DMILP2} outperforms Algorithm  \ref{Alg:DMILP1} with respect to solution time, number of iterations and number of cuts.

\item With the current implementation and numerical results, the main advantage of Algorithms \ref{Alg:DMILP1} and \ref{Alg:DMILP2} is that they preserve data privacy for different blocks. 
\end{enumerate}

Finally, we note that distributed and decentralized optimization are dynamic and evolving area. Data privacy, distributed databases, and computational gains motivate to adapt distributed optimization in many industries such as electric power systems, supply chain, health care systems and etc. Therefore, developing fast and robust distributed exact and heuristic methods for MILPs. 
A possible direction for future research is to blend the speed of R\&F and precision of the exact methods. Another topic for future work is investigating stronger primal cuts to speed up the proposed exact methods. Moreover, the proposed methods can be improved for specific applications by exploiting the problem structures.
 
\bibliographystyle{spmpsci}
\bibliography{references}
\end{document}